\documentclass[12pt]{amsart}
\usepackage{graphicx}

\newcommand{\kt}{K_t}
\newcommand{\kr}{K_{r,f}}

\newcommand{\lcr}{\raisebox{-5pt}{\mbox{}\hspace{1pt}
                 \includegraphics{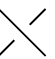}\hspace{1pt}\mbox{}}}
\newcommand{\ift}{\raisebox{-5pt}{\mbox{}\hspace{1pt}
                 \includegraphics{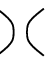}\hspace{1pt}\mbox{}}}
\newcommand{\zer}{\raisebox{-5pt}{\mbox{}\hspace{1pt}
                 \includegraphics{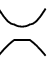}\hspace{1pt}\mbox{}}}

\newcommand{\smst}[1]{\makebox[0pt]{\scriptsize{$#1$}}}

\newtheorem{theorem}{Theorem}
\newtheorem{lemma}{Lemma}
\newtheorem{cor}{Corollary}
\newtheorem{prop}{Proposition}

\title{The Yang-Mills Measure in the Kauffman Bracket Skein Module}

\author{Doug Bullock}

\address{Department of Mathematics, Boise State University, Boise, ID
83725, USA}

\email{\tt bullock@math.idbsu.edu}

\author{Charles Frohman}

\address{Department of Mathematics, University of Iowa, Iowa City, IA
52242, USA}

\email{\tt frohman@math.uiowa.edu}

\author{Joanna Kania-Bartoszy\'{n}ska}

\address{Department of Mathematics, Boise State University, Boise, ID
83725, USA} 

\email{\tt kania@math.idbsu.edu}

\thanks{
This research was partially supported by   NSF-DMS-9803233 and
NSF-DMS-9971905.}

\begin{document}

\begin{abstract} 
For each closed, orientable surface $\Sigma_g$, we construct a
local, diffeomorphism invariant  trace on the Kauffman
bracket skein module $\kt(\Sigma_g \times I)$.  The trace is defined when
$|t|$ is neither $0$ nor $1$, and at certain roots of unity.
At $t = - 1$, the trace
is integration against the symplectic measure on the 
$SU(2)$ character variety of the fundamental group of $\Sigma_g$.
\end{abstract}

\maketitle

\section{Introduction}

Since the introduction of quantum invariants of $3$-manifolds
\cite{RT,Wi1} the fact that they are only defined at roots of unity
has been an obstruction to analyzing their properties. One approach
has been to study the perturbative theory of quantum invariants
\cite{LMO}. However, there is ample evidence quantum invariants
of three manifolds exist as holomorphic functions on the unit disk,
that diverge everywhere on the unit circle but at roots of unity \cite{LZ}.
This paper takes a step towards seeing that this holds in general.
The Yang-Mills measure is the path integral on a
topological quantization \cite{BFK} of the $SU(2)$-characters of the
fundamental group of a closed surface. The measure displays the same
convergence properties as are expected of quantum invariants of
$3$-manifolds.

The Yang-Mills measure in the Kauffman bracket skein algebra of a
cylinder over a closed surface $\Sigma_g$ is a local, diffeomorphism
invariant trace. It quantizes the symplectic measure on the space
$\mathcal{M}(\Sigma_g)$ of conjugacy classes of representations of the
fundamental group of $\Sigma_g$ into $SU(2)$. The definition of the
symplectic structure and formulas for its computation are in
\cite{G1,G2}. The volume of $\mathcal{M}(\Sigma_g)$ was computed by
Witten in \cite{Wi2} in two ways: via the equivalence of two
computations in quantum field theory, and by noting that the
symplectic measure is equal to the measure coming from Reidemeister
torsion. In Witten's setting the
Yang-Mills measure is a path integral in a lattice model of field
theory that depends on area. Forman \cite{F} gave a direct proof that
Witten's measure converges to the symplectic measure as the area goes
to zero.

Alekseev, Grosse and Schomerus \cite{AGS} conceived of a method of
constructing lattice gauge field theory based on a quantum group. This
idea was further developed by Buffenoir and Roche \cite{BuRo} who gave
a construction of the algebra, its Wilson loops and a trace called the
Yang-Mills measure that were completely analogous to Witten's
construction.  Their theory is topological when the area is set to zero.

The method of constructing the algebras in \cite{AGS,BuRo} is
combinatorial and based on generators and relations. 
We gave a new construction of lattice gauge field theory in \cite{la} that is
``coordinate free''. The connections form a co-algebra and the product
on the gauge fields is a convolution with respect to the
co-multiplication of connections. This allows the structure of the
observables to be elucidated. We found working over formal power
series, basing the theory on quantum $sl_2$, that the observables are
the Kauffman bracket skein algebra of a cylinder over a regular
neighborhood of the $1$-skeleton. In \cite{la2} we recover the same
result working over the complex numbers.

These considerations lead one to expect that the Yang-Mills measure
exists as a trace on the Kauffman bracket skein algebra of a closed
surface.  In this paper we affirm this fact, with the only reservation
that if the deformation parameter $t$ is a generic point on the unit
circle, then the measure does not converge.  However, at roots of
unity the trace exists and is well known.
Furthermore, at $t= -1$ the Yang-Mills
measure is the symplectic measure on
$\mathcal{M}(\Sigma_g)$.

This paper is organized as follows. Section 2 recalls definitions,
associated formulas and the algebra structure of the Kauffman bracket
skein module of a cylinder over a surface. In section 3 the Yang-Mills
measure is defined for compact surfaces with boundary, and is proved
to be a trace. In section 4, working with the parameter $t$ such that
$|t|\neq 1$, we obtain estimates for the absolute value of the
tetrahedral coefficients and use these to show that the Yang-Mills
measure can be defined for closed surfaces.  In section 5 we define
and investigate the
measure when $t$ is a root of unity. 

\section{Preliminaries}

Let $M$ be an orientable $3$-manifold.  A framed link in $M$ is an
embedding of a disjoint union of annuli into $M$. Framed links are
depicted by showing the core of an annulus lying parallel to the plane
of the paper (i.e.\ with blackboard framing).
Two framed links in $M$ are equivalent if there is an isotopy of $M$
taking one to the other. Let $\mathcal{L}$ denote the set of
equivalence classes of framed links in $M$, including the empty
link. Fix a complex number $t\neq 0$.  Consider the vector space
$\mathbb{C} \mathcal{L}$ with basis $\mathcal{L}$.  Define $S(M)$ to
be the smallest subspace of $\mathbb{C} \mathcal{L}$ containing all
expressions of the form $\displaystyle{\lcr-t\zer-t^{-1}\ift}$ and
$\bigcirc+t^2+t^{-2}$, where the framed links in each expression are
identical outside balls pictured in the diagrams. The Kauffman bracket
skein module $K_t(M)$ is the quotient
\[ \mathbb{C} \mathcal{L} / S(M). \]
Let $F$ be a compact orientable surface and let $I=[0,1]$. There is an
algebra structure on $K_t(F\times I)$ that comes from laying one link
over the other.  Suppose that $\alpha,\beta \in K_t(F\times I)$ are
skeins represented by links $L_{\alpha}$ and $L_{\beta}$. After
isotopic deformations, to ``raise'' the first link and ``lower'' the
second, $L_{\alpha} \subset F\times (\frac{1}{2},1]$ and $L_{\beta}
\subset F\times [0,\frac{1}{2})$. The skein $\alpha * \beta $ is
represented by $L_{\alpha} \cup L_{\beta}$. This product extends to a
product on $K_t(F\times I)$. We denote the resulting algebra by
$K_t(F)$ to emphasize
that it comes from viewing the underlying three
manifold as a cylinder over $F$.

The notation and the formulas in this paper are taken from \cite{KL}.
However, the variable $t$ replaces $A$, and we use quantum integers
\[ [n] = \frac{t^{2n} - t^{-2n}}{t^2-t^{-2}}. \] When $t=\pm 1$, $[n]=n$.
Note that  $\Delta_n$ from \cite{KL} is equal to $(-1)^n[n+1]$.

There is a standard convention for modeling a skein in $\kt(M)$ on a framed
trivalent graph  $\Gamma \subset M$.  When
$\Gamma$ is represented by a diagram we assume blackboard framing.  An
{\em admissible coloring} of $\Gamma$ is an assignment of a
nonnegative integer to each edge so that the colors at trivalent
vertices form admissible triples (defined below).  The corresponding
skein in $\kt(M)$ is obtained by replacing each edge labeled with the
letter $m$ by the $m$-th Jones--Wenzl
idempotent (see \cite{We}, or \cite{Li}, p.136), and replacing
trivalent vertices with Kauffman triads (see \cite[Fig.\ 14.7]{Li}).

Recall the fusion identity:
\[
\raisebox{-24pt}{\includegraphics{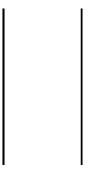}}
\hspace{-20pt}
\raisebox{23pt}{\smst{a}\hspace{23pt}\smst{b}}
= 
\sum_{c} (-1)^c \frac{[c+1]}{\theta(a,b,c)} 
\raisebox{-24pt}{\includegraphics{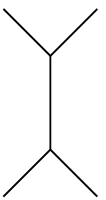}}
\hspace{-9pt} 
\raisebox{3pt}{\smst{c}} 
\hspace{-9pt}
\raisebox{27pt}{\smst{a}\hspace{22pt}\smst{b}}
\hspace{-22pt}
\raisebox{-22pt}{\smst{a}\hspace{22pt}\smst{b}}
\]
where the sum is over all $c$ so that the triples $(a,b,c)$ are
admissible, i.e.  $a+b+c$ is even, $a\leq b+c$, $b\leq a+c$,
and $c\leq a+b$. Value of $\theta(a,b,c)$ is given by equation
(\ref{theta}) below. The fusion relation is satisfied in $\kt(M)$ unless
$t$ is a root of unity other than $\pm 1$.

\section{The Yang-Mills Measure in a Handlebody}

Throughout this section we assume that $t$ is not
a root of unity.
The first result is well known and comes from Przytycki's \cite{P}
construction of examples of torsion in skein modules.

\begin{lemma} 
[{\bf The Sphere Lemma}] Let $s_c$ be a skein represented by coloring
a trivalent framed graph in the manifold $M$.  Suppose further that
there is a sphere embedded in $M$ which intersects the underlying
graph transversely in a single point in the interior of an edge, and
the color of that edge is not zero.  Then $s_c=0$. 
\end{lemma}

\proof Using the ``light bulb trick'' isotope the framed graph $s_c$
so that it is the same graph, but the framing on the edge intersecting
the sphere has been changed by adding two kinks.  Using the formula for 
eliminating a kink,
notice  that $s_c$ is a nontrivial complex multiple of itself. Ergo,
$s_c$ represents zero in $\kt(M)$.  \qed

Consider now $\kt(\#_g S^1\times S^2)$, the Kauffman bracket skein module 
of the connected sum of $g$ copies of $S^1\times S^2$.

\begin{prop}  
The skein module $\kt(\#_g S^1\times S^2)$ is
canonically isomorphic to $\mathbb{C}$. The isomorphism is given by
writing each skein as a complex multiple of the empty skein.
\end{prop}

\proof This follows easily from theorems of Hoste and Przytycki
\cite{HP,P,P2}.  In \cite{HP} the Kauffman bracket skein module of
$S^1\times S^2$ is computed over $\mathbb{Z}[t,t^{-1}]$. This along
with the results in \cite{P2} on the Kauffman bracket skein module of
a connected sum over rational functions in $t$, combined with the
universal coefficient theorem stated in \cite{P}, proves the desired
result.  

We outline the actual isomorphism with the complex numbers.
Choose a system of spheres in $\#_g S^1\times S^2$ that cut it down to
a punctured ball.  
Given a skein in $\#_g S^1\times S^2$, represent it as a linear
combination of colored, framed, trivalent graphs
intersecting the spheres transversely in interior of edges, and so that
each graph intersects any sphere at most once.  This is done by fusing
multiple edges passing through the same sphere. By the sphere lemma,
we can assume the graphs miss the spheres.  Now take the Kauffman
bracket of the skein in the punctured ball to write it as a complex
multiple of the empty skein.  \qed

Given a handlebody $H$ of
genus $g$ its double is $\#_g S^1\times S^2$.  There is a linear
functional $\mathcal{YM}:\kt(H) \rightarrow \mathbb{C}$ computed by
taking the inclusion of $H$ into $\#_g S^1\times S^2$ followed by
taking the `` Kauffman bracket'' as above.  Let $F$ be a compact,
oriented surface with boundary. Since $F\times I$ is a handlebody the
linear functional
\[ \mathcal{YM}:\kt(F) \rightarrow \mathbb{C},\]
is defined. We  call this the {\em Yang-Mills measure}.

Choose a trivalent spine of $F$.  The admissible colorings of that
spine form a basis for $\kt(F)$. The skein modules of the disk 
and annulus are exceptions; the first is spanned by
the empty skein and the latter is described in Section 4.  In terms
of this basis the Yang-Mills measure is just the coefficient of the
skein coming from labeling all the edges of the spine with $0$.

\begin{prop} 
The Yang-Mills measure is a trace, that is
\[\mathcal{YM}(\alpha *\beta)=\mathcal{YM}(\beta *\alpha).\]
Furthermore, the trace is invariant under the action of the
diffeomorphisms of $F\times I$ on $\kt(F)$.
\end{prop}

\proof Let $L$ be the link $\partial F \times \{1/2\}$. The result of
removing $L$ from the double of $F \times I$ is homeomorphic to the
Cartesian product of the interior of $F$ with a circle.  Given any
skein in $F\times I$ we can represent it by a linear combination of
framed links that miss $L$.  Hence, the Yang-Mills measure factors
through the skein module of $F \times S^1$. In $F\times S^1$ the
skeins $\alpha *\beta$ and $\beta * \alpha$ are the same.

The group of diffeomorphisms of the handlebody $F\times I$ acts on
$\kt(F)$ in the obvious way. If $f: F\times I \rightarrow F\times I$
is a diffeomorphism then it can be extended to $Df:\#_g S^1\times S^2
\rightarrow \#_g S^1 \times S^2$.  Since the image of the empty skein
under a diffeomorphism is the empty skein, the action of $Df$ on
$\kt(\#_g S^1\times S^2)$ is trivial.  Therefore,
$\mathcal{YM}(f(\alpha))=\mathcal{YM}(\alpha)$. \qed

The final commonly used property of the Yang-Mills measure is that it
is {\em local}.  Suppose that $k$ is a proper arc in $F$. Cut $F$
along $k$ to get a surface $F'$.  It is evident that if we write a
skein $\alpha$ as a linear combination of admissibly colored graphs,
each one intersecting $k$ transversely in at most a single point, then we can
throw out any graph such that the edge intersecting $k$ carries a
nonzero label. This yields a skein in $F'$, denoted by
$\alpha_k$. Then $\mathcal{YM}(\alpha)=\mathcal{YM}(\alpha_k)$.

\section{The Yang-Mills measure on a closed surface}

Throughout this section assume that $|t|\neq 1$. In fact, we only work
with $0<t<1$. However, it is evident that the same proofs are valid
when $1<t$ since the formulas are symmetric in $t$ and
$t^{-1}$. Finally, the arguments extend to the case where $t$ is not
real by replacing the estimates for $t\in {\mathbb R}$ by estimates of
the absolute value of $t\in {\mathbb C}$.

Recall the Kauffman bracket skein algebra of a cylinder over an
annulus $A$. The central core of the annulus can be seen as a link by
giving it the blackboard framing.  Let $s_i$ be the skein in the
annulus which is the result of plugging the $i$-th Jones-Wenzl
idempotent into the core. The skein module $K_t(A)$ is the vector
space with basis $\{s_i\}$, where $i$ runs from zero to infinity. The
product with respect to this basis is given by
\begin{equation}\label{mult}
s_i * s_j=\sum_{q\geq |i-j|, \text{by 2's}}^{i+j}s_q.
\end{equation}

Use the Yang-Mills measure on $\kt(A)$ to define a pairing:
\begin{equation}\label{pair}
\langle\alpha,\beta\rangle = \mathcal{YM}(\alpha * \beta).
\end{equation}
The $s_i$ form an orthonormal basis with respect
to (\ref{pair}). This pairing identifies the  linear dual of 
$K_t(A)$ with
series of the form  $\sum_i \alpha_i s_i$, where the $\alpha_i$ are complex
numbers. Note that:
\[\langle\sum_{i=0}^{\infty} \alpha_i s_i,\sum_{j=0}^n \beta_j
s_j\rangle
=\sum_{i=0}^n \alpha_i\beta_i.\]
Let $\Sigma_{g,1}$ denote the compact orientable surface of genus $g$
with one boundary component. There is a pairing,
\[ K_t(A) \otimes K_t(\Sigma_{g,1}) \rightarrow  K_t(\Sigma_{g,1})\]
given by representing the skein in $K_t(\Sigma_{g,1})$ by a linear
combination of links disjoint from some collar of the boundary, and
plugging the skein in $K_t(A)$ into the collar. The Yang-Mills
measure can then be applied to give a pairing,
\begin{equation}\label{pairing}  
K_t(A) \otimes K_t(\Sigma_{g,1}) \rightarrow  \mathbb{C}.
\end{equation}
This means there is a well defined map,
\[ Y:K_t(\Sigma_{g,1})\rightarrow K_t(A)^*.\]
Topologize $K_t(A)$ by giving it the weak topology from $Y$. That is a
sequence  $\sigma_n \in K_t(A)$ is Cauchy if for every skein $\alpha
\in K_t(\Sigma_{g,1})$, the sequence of complex numbers
$Y(\alpha)(\sigma_n)$ is Cauchy.  A linear functional on
$K_t(\Sigma_{g,1})$ that comes from an element of this completion via
the pairing (\ref{pairing}) is
called a {\em distribution}. It is interesting to note that the weak
topology from $Y$ on $\kt(A)$ depends on the genus of the surface.

If $g>1$ there is a distribution on $K_t(\Sigma_{g,1})$ which
annihilates all ``handle-slides'' ({\em Skeins that are represented by
the difference of two links such that one can be obtained from the
other by a slide 
across an imagined disk filling the boundary of $\Sigma_{g,1}$}).
This linear functional descends to the skein module of the closed
surface.  Yang-Mills measure on a closed surface is the result of
evaluating this distribution followed by a normalization.

Let's think about what a skein in $\kt(A)$ would be like if it
annihilated all handle-slides. Begin by writing it as $\sum_i \alpha_i
s_i$ and solve for the $\alpha_i$. A simple computation shows that if
$\alpha_0$ is  zero then all $\alpha_i$ are zero. Normalize so that 
$\alpha_0=1$.  Notice
that if our skein annihilates handle-slides then the skein
$s_1+[2]s_0$ must be annihilated. Using the rules for multiplication (\ref{mult})
we see that the coefficient $\alpha_1$ is equal to $-[2]$. Continuing
on this way we see that this skein has to be \[ \sum_i (-1)^i[i+1]s_i,\] which
is of course not in $K_t(A)$.

The first goal is to show that for $g>1$ the sequence of partials sums
$\sum_{i=0}^n (-1)^i[i+1]s_i$ is Cauchy in the weak topology from $Y$,
and so defines a distribution.

\begin{figure}
\begin{picture}(56,120)
\raisebox{8pt}{\hspace{-1in}\includegraphics{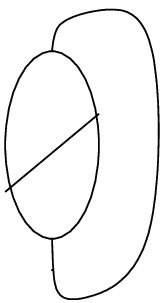}
\put(-50,16){$a$}
\put(-60,49){$b$}
\put(5,49){$e$}
\put(-23,70){$c$}
\put(-15,40){$d$}
\put(-40,49){$f$}}
\end{picture}
\begin{picture}(46,45)
\hspace{0.1in}\raisebox{40pt}{\includegraphics{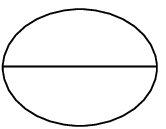}
\put(-25,5){$a$}
\put(-25,19){$b$}
\put(-25,37){$c$}}
\end{picture}
\caption{Tet and theta}\label{tet}
\end{figure}
The notation $\text{Tet}\begin{pmatrix} a & b & e \\ c & d &
f\end{pmatrix}$ stands for the Kauffman bracket of the skein pictured
in Figure \ref{tet} on the left.  The explicit formula is given in
\cite{KL}. We also need the quantity $\theta(a,b,c)$ which is the
Kauffman bracket of the colored graph on the right in Figure
\ref{tet}.  In terms of quantum integers
\begin{equation}\label{theta}
\theta(a,b,c)= (-1)^{\frac{a+b+c}{2}} \frac{[\frac{a+b+c}{2}+1]!
  [\frac{a+b-c}{2}]! [\frac{b+c-a}{2}]! [\frac{c+a-b}{2}]!}
{[a]![b]![c]!}.
\end{equation}

Another quantity, called a $6j$ symbol, is derived from the
tetrahedral evaluation. Specifically,
\begin{equation}\label{6j} 
\left\{ \begin{matrix} a & b & e \\ c & d & f \end{matrix}\right\} =
\frac{\text{Tet}\begin{pmatrix} a & b & e \\ c & d &
f\end{pmatrix}(-1)^e[e+1]} {\theta(a,d,e)\theta(c,b,e)}.
\end{equation}
The $6j$
symbols can be woven together to give a change of basis matrix for the
Whitehead move on graphs. As a consequence they satisfy an
orthogonality equation:

\begin{equation}\label{ort}
 \sum_e  \left\{ \begin{matrix} a & b & e \\ c & d & f \end{matrix}\right\}
 \left\{ \begin{matrix} d & a & g \\ b & c & e
   \end{matrix}\right\}=\delta_f^g,
\end{equation}
where $\delta_f^g$ is the Kronecker delta.

The following proposition seems quite weak, but turns out to be a powerful
tool for gauging the convergence of series of Kauffman brackets.

\begin{prop}\label{est1}
\[ \left| \text{Tet}\begin{pmatrix} a & b & e \\ c & d & f\end{pmatrix}\right| 
\leq
\sqrt{\frac{\theta(b,c,e)\theta(a,d,e)\theta(a,b,f)\theta(c,d,f)}{(-1)^{e+f}[e+1][f+1]}}\]
\end{prop}

\proof In order for all the triples at the vertices of a tetrahedron
to be admissible , the parity of the sum of the entries in any two
columns of
\[ \text{Tet}\begin{pmatrix} a & b & e \\ c & d & f\end{pmatrix}\]
has to be the same. Use (\ref{6j}) to expand the formulas for the $6j$ symbols in the
orthogonality relation (\ref{ort}), with
$g=f$. The tetrahedral evaluations are equal and the signs of the
$\theta$'s and the $(-1)^{e+f}$ cancel so that each term in the sum is
positive.  Hence every term in the sum is less than $1$. Fixing $e$ and
putting everything except for the tetrahedral evaluations on the right
hand side, and taking square roots yields the desired result.\qed

\begin{cor}\label{est2} 

There is a real valued function $C(k_1,k_2,k_3)$ so that
\begin{equation}\label{term}
\frac{|\text{Tet}\begin{pmatrix} i & i & i \\ k_1 & k_2 & k_3
\end{pmatrix}|}
{\sqrt{|\theta(i,i,k_1)\theta(i,i,k_2)\theta(i,i,k_3)}|}
\end{equation}
is less than $t^{i}C(k_1,k_2,k_3)$ whenever the graphs corresponding
to the functions in the formula are admissibly labeled.\end{cor}

\proof Substitute into the inequality from Proposition \ref{est1} to
get,
\begin{equation}\label{fund}
\left| \text{Tet}\begin{pmatrix} i & i & i \\ k_1 & k_2 &
k_3\end{pmatrix}\right| \leq \sqrt{\frac{\theta(k_1,k_2,k_3)
\theta(i,i,k_1) \theta(i,i,k_2)
\theta(i,i,k_3)}{(-1)^{i+k_3}[k_3+1][i+1]}}.
\end{equation} 
Shift
$\sqrt{\theta(i,i,k_1)\theta(i,i,k_2)\theta(i,i,k_3)}$ to the left
hand side.  Use the fact that $\frac{1}{[i+1]} \leq  t^{2i}$ to make the right
hand side bigger. Finally, note that the remaining factor on the right hand
side is a function of $k_1$, $k_2$ and $k_3$. \qed

\begin{theorem} The sequence $ \sum_{i=0}^n
(-1)^i[i+1]s_i$ defines a distribution for $g>1$. That is, the limit
\[\mathcal{YM}_D(\alpha)= 
\lim_{n\rightarrow \infty} \mathcal{YM}(\alpha * \sum_{i=0}^n
(-1)^i[i+1]s_i)\] exists and gives a  well defined trace
 on $\kt(\Sigma_{g,1})$ when
$g>1$. \end{theorem}

\proof Choose a trivalent spine for $\Sigma_{g,1}$ with $4g-2$
vertices and $6g-3$ edges. Basis elements $s_c$ for
$\kt(\Sigma_{g,1})$ correspond to labeling the edges admissibly with
integers $k_j$, where $j$ runs from $1$ to $6g-3$. Let $s_i$ denote
the core of an annulus that runs parallel to the boundary, labeled
with the $i$th Jones-Wenzl idempotent. In order to compute
$\mathcal{YM}(s_c*s_i)$ place  both skeins in the same
diagram. Choose a system of arcs, each intersecting this configuration
transversely in three points, that isolate the vertices from one
another. The transverse points of intersection are labeled $i$, $k_j$,
$i$ as you traverse each arc. Fuse along these arcs, until
the resulting graphs intersect each arc in at most one point. Discard any
term where the label on an edge intersecting an arc is not zero. Given a vertex $v$, let
$(k_{v1},k_{v2},k_{v3})$ be the triple of colors appearing there. The
resulting answer is:
\begin{equation}\label{prod} \mathcal{YM}(s_c*s_i)=
\prod_{j=1}^{6g-3} \frac{1} {\theta(i,i,k_j)} \prod_{v} \text{Tet}
\begin{pmatrix} i & i & i \\ k_{v1} & k_{v2} & k_{v3} \end{pmatrix}.
\end{equation}

Each edge appears at exactly two vertices, so (\ref{prod}) can be written
as a product of  $4g-2$ factors like (\ref{term}). By Corollary
\ref{est2} the absolute value of $\mathcal{YM}(s_c*s_i)$ is less than
$C(k_j)t^{i(4g-2)}$, where $C(k_j)$ is a number depending only on the
$k_j$. The $n$th partial sum for $\mathcal{YM}_D(s_c)$ is
\[ 
\sum_{i=0}^n (-1)^i[i+1]\prod_{j=1}^{6g-3} \frac{1} {\theta(i,i,k_j)} 
\prod_{v} \text{Tet}
\begin{pmatrix} i & i & i \\ k_{v1} & k_{v2} & k_{v3} \end{pmatrix}.
\]
Note that $[i+1]$ is less than
$(i+1)t^{-2i}$. Hence the $i$-th summand is less than
$(i+1)(-1)^iC(k_j)t^{i(4g-4)}$. The ratio test implies that the
sequence of partial sums is absolutely convergent for $0<t<1$. 

Finally,
$\mathcal{YM_D}$ is a trace since the partial sums
$\sum_{i=0}^n(-1)^i[i+1]s_i$ can be seen as lying in the center of
$\kt(\Sigma_{g,1})$. \qed

\begin{theorem} 
$\mathcal{YM}_D$ descends to give a well defined trace
\[ \mathcal{YM}: \kt(\Sigma_{g})\rightarrow \mathbb{C}.\]
\end{theorem}

\proof There is a homomorphism $K_t(\Sigma_{g,1}) \rightarrow
\kt(\Sigma_{g})$ induced by inclusion.  The 
surface $\kt(\Sigma_{g})$ is the result of adding a disk to the
boundary of surface 
$K_t(\Sigma_{g,1})$. The kernel of this homomorphism consists of all skeins that
can be written as a linear combination of handle-slides.
The next step is to show
that the linear functional $\mathcal{YM}_D$ annihilates all
handle-slides. To this end we analyze the difference of the two skeins
in the annulus (relative to a pair of points in the boundary).
\begin{equation}\label{annulus}
 \sum_{i=0}^n(-1)^i[i+1]
\left(\raisebox{-.25in}{\includegraphics{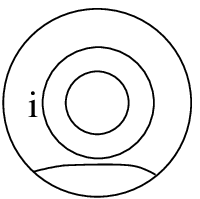}}-
\raisebox{-.25in}{\includegraphics{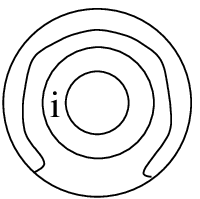}} \right)
\end{equation} 

The analysis of the diagram (\ref{annulus}) diagram is due to
Lickorish, \cite{Li}. It is equal to:
\begin{equation}\label{nosum}
(-1)^n[n+1] \left(\raisebox{-.25in}{\includegraphics{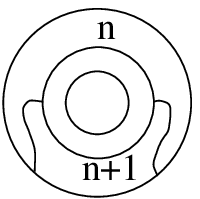}}-
\raisebox{-.25in}{\includegraphics{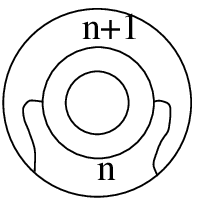}} \right).
 \end{equation}

This diagram needs to be set in place. Using standard arguments as in
\cite{Bu} yields that we only need to check handle-slides of the
following form. Take a skein corresponding to a colored spine, and
separate one strand along an edge.

\begin{picture}(219,85)
\hspace{1.87in}\scalebox{.75}{\includegraphics{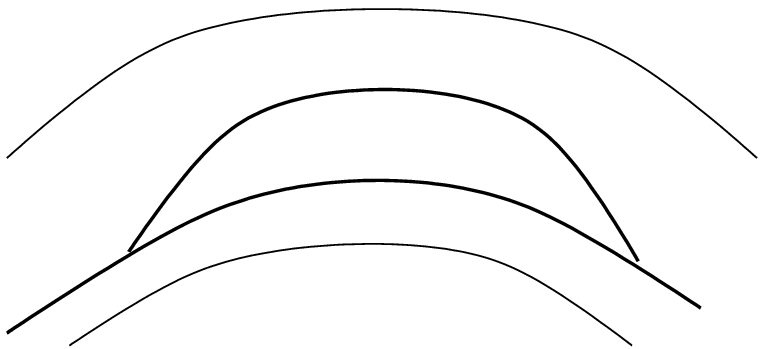}
\put(-200,25){$k$} \put(-120,50){$k-1$} \put(-20,25){$k$}}
\end{picture}

Now slide the strand over the added disk, locally the diagram looks
like:

\begin{picture}(219,95)
\hspace{1.87in}\scalebox{.75}{\includegraphics{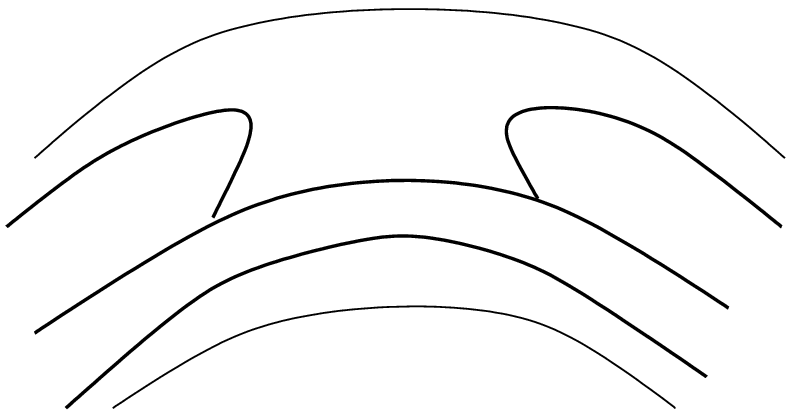}
\put(-200,45){$k$}
\put(-120,70){$k-1$}
\put(-20,40){$k$}}
\end{picture}

Multiplying each of the diagrams above by $\sum_{i=0}^n (-1)^i[i+1]
s_i$, taking their difference, and using the identity
(\ref{annulus})=(\ref{nosum}), we get a difference of two terms like
the one below. In the first one the label $u=n$ and the label $v=n+1$,
and in the second one $u=n+1$ and $v=n$.

\begin{picture}(219,85)
\hspace{1.87in}\scalebox{.75}{\includegraphics{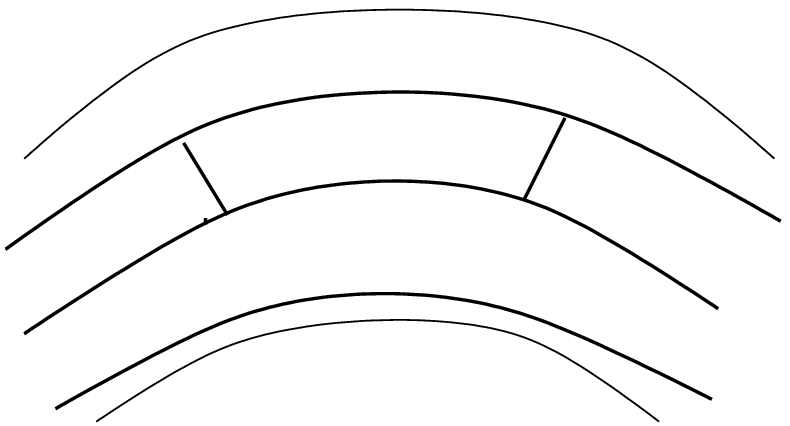}
\put(-200,25){$k$} \put(-120,50){$k-1$} \put(-120,100){$v$}
\put(-200,74){$u$} \put(-20,70){$u$} \put(-20,25){$k$}
\put(-20,10){$u$}}
\end{picture}

Fusing to isolate the vertices of this diagram requires two more cross
cuts than the diagrams we have been working with up till now. We get
the product of
\begin{equation}\label{newpart}
(-1)^n [n+1] \frac{1} {\theta(u,k,u) \theta(u,k-1,v) } \text{Tet}
\begin{pmatrix} u & u & v \\  1 & k-1 & k \end{pmatrix}\text{Tet}
\begin{pmatrix} u & v & u \\   1 & k & k-1 \end{pmatrix}  
\end{equation}   
with the standard product,
\begin{equation}\label{oldpart}
\prod_{j=1}^{6g-3}\frac{1}{\theta(u,u,k_j)}\prod_{v}\text{Tet}\begin{pmatrix}
u & u & u \\ k_{v1} & k_{v2} & k_{v3} \end{pmatrix}.
\end{equation}
The product (\ref{oldpart}) is smaller
than a global constant, depending on the $k_j$, times $t^{n(4g-2)}$. It
remains to ascertain that the term (\ref{newpart}) is not too large.  Using
the inequality from Proposition \ref{est1} we get that, regardless of
whether $u=n$ and $u=n+1$, or $u=n+1$ and $u=n$, the absolute value of
(\ref{newpart}) is less than $[n+2]$, which is a universal constant
times $t^{-2n}$. As long as the genus of the surface is greater than
$1$, the full product goes to zero as $n$ goes to infinity. So, in the
limit, all handle-slides are annihilated.  \qed

The case of a surface of genus $1$ is slightly different. To get a
convergent distribution we need to divide the partial sum
$\sum_{i=0}^n (-1)^i[i+1] s_i$ by $n$. The sequence is then Cauchy and
defines a distribution on $\kt(T^2)$.  

The algebra $K_t(T^2)$ is very
nice for working examples. If $(p,q)$ is a pair of integers that are
relatively prime there is an obvious skein $s_{(p,q)}$ in $K_t(T^2)$
corresponding to the $(p,q)$ curve on the torus .  Define a family of
skeins based on $(p,q)$ by using the following iterative scheme:
$s_{(p,q)_0}=2s_{(0,0)}$, that is, twice the empty skein, and
$s_{(p,q)_1}=s_{(p,q)}$. For $d>1$ define:
\[ s_{(p,q)_d}=s_{(p,q)}* s_{(p,q)_{d-1}}- s_{(p,q)_{d-2}}.\]
Finally, if $d=\text{gcd}\{p,q\}$, let \[s_{(p,q)}=s_{(p/d,q/d)_d}.\]
Using this notation the product in $K_t(T^2)$ is given by
\begin{equation}\label{gelca}
s_{(p,q)}*s_{(u,v)}=t^{\left|\begin{matrix} p & q \\ u &
v\end{matrix}\right|} s_{(p+u,q+v)} +t^{-\left|\begin{matrix} p & q \\
u & v\end{matrix}\right|}s_{(p-u,q-v)}.
\end{equation}
The formula (\ref{gelca}) is proven in \cite{FG}.

There is a map
\[\mu : K_t(T^2) \rightarrow \mathbb{C}\emptyset \oplus \mathbb{C}H_1(T^2;Z_2)\] 
introduced in \cite{mike}. Let
\[\mu\left( \sum_{(p,q)} a_{(p,q)}s_{(p,q)}\right)=a_{(0,0)}\emptyset +
\sum_{(p,q)\neq (0,0)} a_{(p,q)}[(p,q)],\] where $[(p,q)]$ is the
$Z_2$--homology class in $H_1(T^2;Z_2)$ corresponding to
$d=\text{gcd}\{p,q\}$ copies of a $(p/d,q/d)$ curve on the torus.  The
map $\mu$ has as its kernel the submodule of all commutators.  Hence
any linear functional on the five dimensional space that is the image
of $\mu$ is a trace. It is easy to check that there is a three
dimensional family of traces 
that are  invariant under diffeomorphism. In this set up
\[
\mathcal{YM}\left(\sum_{(p,q)} a_{(p,q)}s_{(p,q)}\right)=a_{(0,0)}. 
\]
This is the same trace as the one induced from the
inclusion of $\kt(T^2)$ into the non-commutative torus \cite{FG}.
 
Towards uniqueness of the Yang-Mills measure, it should be normalized,
just as the symplectic measure on moduli space needs to be normalized.
It should also be invariant under diffeomorphism, and be local.
Locality is made up by two rules. One for cutting a surface along an
arc and one for removing a point from a closed surface.  If we
formalize our rules correctly, we get the following:
\begin{theorem} 
The Yang-Mills measure is the unique, local, diffeomorphism invariant
trace on $\kt(\Sigma_g)$ up to normalization.\qed
\end{theorem}

\section{Roots of Unity}

Fusion no longer holds in $\kt(M)$ when $t$ is a root of
unity. However, when $t=e^{\frac{\pi i}{2r}}$ then one can take a
quotient, where an appropriate form of the fusion identity is
true. This can be done by setting any skein containing the $(r-1)$-st
Jones-Wenzl idempotent equal to zero.  The quotient is denoted
$\kr(M)$. The {\em reduced skein} is a central object in the
construction of quantum invariants of $3$-manifolds
\cite{FK,Ro1,Ro2}. 

The Yang-Mills measure on a surface with boundary
is obtained the same way as for other values of $t$.
Since $[r]=0$, the iterative procedure for finding a skein in the
annulus that annihilates handle-slides terminates, to yield
\[ \sum_{i=0}^{r-2} (-1)^i[i+1] \bigcirc^i.\]
There is  an induced trace,
\[ \mathcal{YM}:\kr(\Sigma_g) \rightarrow \mathbb{C},\]
constructed the same way as for other $t$ except that there is no need
to take a limit because the formula is a finite sum.

Notice that $\Sigma_g$ is the boundary of a handlebody $H_g$ (it
doesn't make any difference which one).  There is an action of
$\kr(\Sigma_g)$ on $\kr(H_g)$ given by gluing skeins in $\Sigma_g
\times I$ into a collar of the boundary of $H_g$.  The action gives a
map
\[ \phi: \kr(\Sigma_g) \rightarrow \mathrm{End}(\kr(H_g)).\]
As we are working at a root of unity, $\kr(H_g)$ is a finite
dimensional vector space. Denote its dimension by $d$, and let
$\omega=\mathcal{YM}(\emptyset)=\sum_{i=0}^{r-2}\frac{1}{[i+1]^{2g-2}}$. 
The Yang-Mills measure is:
\[ \mathcal{YM}(\alpha)= \frac{\omega}{d}\mathrm{tr}(\phi(\alpha)).\]

From \cite{Ro3} the map $\phi$ is injective and
onto. Hence we can identify $\kr(\Sigma_g)$ with
$\mathrm{End}(K_{r,f}(H_g))$.  The Yang-Mills measure is zero on
commutators. Thus it factors through
\[\mathrm{End}(K_{r,f}(H_g)) /[\mathrm{End}(K_{r,f}(H_g)),
\mathrm{End}(K_{r,f}(H_g))].\] 
This quotient is 
a $1$-dimensional vector space. Hence any two linear functionals that
factor through this quotient are equal if they agree on the identity
matrix. The trace also vanishes on commutators, thus it  factors
through the commutator quotient. The normalization in the formula
causes the two induced linear functionals to be the same.

Next we address the cases of $t=\pm 1$. Since the formula for the measure of
a spine is an even function of $t$,  we only need to consider one
value. The value $t=-1$ is more convenient as the correspondence
between $K_{-1}(F)$ and the $SU(2)$-characters of $\pi_1(F)$ is
simpler. The skein of the disjoint union of curves $c_i$ corresponds to
the function that sends the representation $\rho$ to
\[ \prod_i -\mathrm{tr}(\rho(c_i)).\]

\begin{theorem} The Yang-Mills measure is well defined on $K_{\pm1}(\Sigma_g)$
for $g>1$. Let
$s_c$ be an admissibly colored trivalent
spine of $\Sigma_g$.
If $t_n$, with $|t_n|\neq 1$, is a sequence of
complex numbers converging to $\pm 1$ then
\[ \lim_{n\rightarrow \infty} \mathcal{YM}_{t_n}(s_c)  =\mathcal{YM}_{\pm 1}(s_c).\]
\end{theorem}

\proof The formulas for working with skeins in $K_{-1}(F)$ are the
same as the ones for $|t|\neq 1$ except that  quantized integers are replaced
by ordinary integers.  These formulas are the limits as $t\rightarrow
-1$ of the values we have been using. Revisiting the fundamental
estimate (\ref{fund}), we see that,
\begin{equation}
\frac{|\text{Tet}\begin{pmatrix} i & i & i \\ k_1 & k_2 & k_3
\end{pmatrix}|}
{\sqrt{|\theta(i,i,k_1)\theta(i,i,k_2)\theta(i,i,k_3)}|}
\leq
\sqrt{\frac{\theta(k_1,k_2,k_3)}{(-1)^{i+k_3}(k_3+1)(i+1)}}
\end{equation}
from which we conclude that the right hand side is less than or equal to
\[ \frac{C(k_1,k_2,k_3)}{\sqrt{i+1}}.\]
Considering the series for the Yang-Mills measure of a spine, 
comparison to the p-series  implies that it converges as long as the surface has genus
greater than $1$. Similarly, the Yang-Mills measure is invariant under
handle-slides.

The convergence statement follows from the fact that the series that
define the Yang-Mills measure at $t_n$ converge absolutely, and the
terms of the series converge to the terms of the series for the
Yang-Mills measure at $-1$. \qed

For a surface of genus $1$ we divide the partial sums, as before, by
the number of terms in the sum, and the series then converges.

\begin{theorem} 
The Yang-Mills measure at $t=-1$
is the symplectic measure on $\mathcal{M}(\Sigma_g)$.
\end{theorem}

\proof Using Weyl orthogonality to compute Witten's Yang-Mills measure
for a surface of area $\rho$ yields that its value on the spine $s_c$
is given by the series
\[\sum_{i=0}^\infty (-1)^i(i+1)e^{-\rho c_2(i)}\prod_{j=1}^{6g-3} 
\frac{1}{\theta(i,i,k_j)} \prod_{v} \text{Tet}
\begin{pmatrix} i & i & i \\ k_{v1} & k_{v2} & k_{v3} \end{pmatrix},\]
where the edges of $s_c$ carry colors $k_i$, and $k_{v_i}$ are the
colors of the edges ending at the vertex $v$, and $c_2(i)$ is the
value of the quadratic Casimir operator on the $(i+1)$-dimensional
irreducible representation 
of $SU(2)$. As both Witten's series and our series
converge absolutely, and Witten's formula converges term by term to
our formula as $\rho \rightarrow 0$, the limit of Witten's Yang-Mills
measure is equal to our Yang-Mills measure at $t=-1$. Finally, Forman
\cite{F} showed that the limit as $\rho \rightarrow 0$ of Witten's
measure is the symplectic measure on $\mathcal{M}(\Sigma_g)$,
normalized as in \cite{F}. \qed

Suppose now that $|t|=1$ and $t$ is not a root of unity. Evaluation of
the Yang-Mills measure on the empty skein on a surface of genus $g$
yields  $\sum_{i=o}^{\infty}\frac{1}{[i+1]^{2g-2}}$. As $t$
is not a root of unity the number $[i+1]^{2g-2}$ gets arbitrarily
close to $1$ infinitely often, which means that the series does not
converge. Therefore the Yang-Mills measure does not exist away from
roots of unity on the unit circle.


\begin{thebibliography}{9999}


\bibitem{AGS} A.Y. Alekseev, H. Grosse, V. Schomerus, {\em
Combinatorial Quantization of the Hamiltonian Chern-Simons Theory I,II},
 Comm. Math. Phys.{\bf 172} (1995), no. 2, 317--358, and
 Comm. Math. Phys.{\bf 174} (1996), no. 3, 561--604.



\bibitem{Bu} D. Bullock, {\em The $(2,\infty)$-skein module of the
complement of a $(2,2p+1)$ torus knot.}, J. Knot Theory Ramifications,
{\bf 4} (1995), no. 4, 619--632.

\bibitem{BFK} D. Bullock, C. Frohman, J. Kania-Bartoszy\'{n}ska, {\em
Understanding the Kauffman bracket skein module}, JKTR, {\bf 8}
(1999), 265--277.

\bibitem{la} D. Bullock, C. Frohman, J. Kania-Bartoszy\'nska, {\em
Topological interpretations of Lattice Gauge Field Theory}, Comm.\
Math.\ Phys. {\bf 198} (1998) 47--81.

\bibitem{la2} D. Bullock, C. Frohman, J. Kania-Bartoszy\'nska, {\em
The Kauffman Bracket Skein as an Algebra of Observables}, preprint.

\bibitem{BuRo} E. Buffenoir, Ph. Roche, {\em Two Dimensional Lattice
Gauge Field Theory Based on a Quantum Group}, Comm. Math. Phys. {\bf
170} (1995), 669-698.

\bibitem{F} R. Forman,  {\em Small volume limits of $2$-d
Yang-Mills} Comm. Math. Phys. {\bf 151} (1993), no. 1, 39--52.

\bibitem{FG} C. Frohman, R. Gelca, {\em Skein modules and the
noncommutative torus}, Transactions of the AMS, to appear.

\bibitem{FK} C. Frohman, J. Kania-Bartoszy\'{n}ska, {\em A Quantum
Obstruction to Embedding}, Math. Proc. of the Cambridge Philosophical
Society, to appear.

\bibitem{G1} W. M. Goldman, {\em The symplectic nature of fundamental
groups of surfaces} Adv. in Math., Advances in Mathematics {\bf 54}
(1984), no. 2, 200--225.

\bibitem{G2} W. M. Goldman, {\em Invariant functions on Lie groups and
Hamiltonian flows of surface group representations}, Invent. Math.{\bf
85} (1986), no. 2, 263--302.

\bibitem{HP} J. Hoste, J. Przytycki, {\em The Kauffman bracket skein
module of $S^1\times S^2$}, Mathematische Zeitschrif, {\bf 220}
(1995), no. 1, 65--73.
 
\bibitem{KL} L. H. Kauffman and S. Lins, {\em Temperley-Lieb
recoupling theory and invariants of $3$-manifolds}, Ann.\ of Math.\
Studies {\bf 143}, Princeton University Press (1994).

\bibitem{LZ} R. Lawrence and D.  Zagier, {\em Modular forms and
quantum invariants of $3$-manifolds}, in: Sir Michael Atiyah: a great
mathematician of the twentieth century.  Asian J. Math. {\bf 3},
(1999) 93--107.

\bibitem{Li} W. B. R. Lickorish, {\em An Introduction to Knot Theory},
Springer, GTM {\bf 175}, 1997.

\bibitem{mike} M. McLendon, personal communication.

\bibitem{LMO} T. Le, H. Murakami, J. Murakami, T. Ohtsuki, {\em A
three-manifold invariant via the Kontsevich integral}, Osaka
J. Math. 36 (1999), no. 2, 365--39

\bibitem{P} J. H. Przytycki, {\em Fundamentals of Kauffman Bracket
Skein Modules}, httt://xxx.lanl.gov/math.GT/9809113.

\bibitem{P2} J. H. Przytycki, {\em Kauffman bracket skein module of a
connected sum of 3-manifolds} httt://xxx.lanl.gov/math.GT/9911120.

\bibitem{RT} N. Y. Reshetikhin, V. G. Turaev, {\em Invariants of
$3$-manifolds via link polynomials and quantum groups},
Invent. Math. {\bf 103}, (1991) 547--597.

\bibitem{Ro1} J. Roberts, {\em Skein theories as TQFTs}, preprint.

\bibitem{Ro2} J. Roberts, {\em Skein theory and Turaev-Viro
invariants}, Topology {\bf 34} (1995) 771--787.

\bibitem{Ro3} J. Roberts, {\em Quantum Invariants via Skein Theory},
Thesis, Pembroke College, Cambridge, 1994.

\bibitem{We} H. Wenzl, {\em On sequences of projectors},
  C.R. Math. Rep. Acad. Sci. {\bf IX}, (1987) 5-9.

\bibitem{Wi1} E. Witten, {\em Quantum field theory and the Jones
polynomial}, Comm. Math. Phys. {\bf 121}, (1989) 351-399.

\bibitem{Wi2} E. Witten, {\em Quantum Gauge Theories in Dimension
Two}, Comm. Math. Phys. {\bf 141}, (1991) 153-209.

\end{thebibliography}
\end{document}